\documentclass[12pt,a4paper,oneside]{amsart}
\usepackage[left=1.1in, right=0.8in, top=1.0in, bottom=1.0in]{geometry}
\usepackage{graphicx}
\usepackage{setspace}
\usepackage{indentfirst} 
\usepackage{cases}
\usepackage{wrapfig}
\usepackage[toc,page]{appendix}
\usepackage{amsmath}
\usepackage{amsthm}
\usepackage{amssymb}
\usepackage{enumerate}
\usepackage{mathrsfs}
\usepackage{dcolumn}
\usepackage{xcolor}
\usepackage{mathtools}
\usepackage{bookmark}
\usepackage[nocompress]{cite}

\newcolumntype{L}{D{.}{.}{2,5}}
\theoremstyle{plain}
\newtheorem{thm}{Theorem}[section]
\newtheorem*{mainthm}{Milman--Pettis Theorem}
\newtheorem{lemma}{Lemma}[section]
\newtheorem{remark}{Remark}[section]

\newtheorem{example}{Example}[section]
\newtheorem{corollary}{Corollary}[section]
\newtheorem{proposition}{Proposition}[section]

\DeclareMathOperator{\cl}{cl}

\theoremstyle{definition}
 \usepackage{hyperref}
 \hypersetup{
    colorlinks=true,
    linkcolor=black,
    filecolor=magenta,      
    urlcolor=cyan}
  \usepackage[capitalize,nameinlink]{cleveref}

\begin{document}
	\title[Banach spaces of sequences arising from infinite matrices]{Banach spaces of sequences arising \\ from infinite matrices}
	\author{{A. B\"erd\"ellima}$^*$ and {N. L. Braha}$^{\dagger}$}
	\thanks{$^*$ berdellima@gmail.com, $\dagger$ nbraha@gmail.com}
	\AtEndDocument{\bigskip{%
	  \textsc[$^*$]{German International University in Berlin,\\ Faculty of Engineering, Berlin 13507, Germany.} \par
	  \addvspace{\medskipamount}
	  \textsc[$\dagger$]{University of Prishtina,\\ Faculty of Natural Science and Mathematics, \\Department of Mathematics and Computer Sciences,\\ Av. Mother Teresa, Nr.5, 10000 Prishtina, Kosovo.} \par	}}

%
	\maketitle
	
	\begin{abstract}
	Given an infinite matrix $M=(m_{nk})$ we study a family of sequence spaces $\ell_M^p$ associated with it. When equipped with a suitable norm $\|\cdot\|_{M,p}$ we prove some basic properties of the Banach spaces of sequences $(\ell_M^p,\|\cdot\|_{M,p})$. In particular we show that such spaces are separable and strictly/uniformly convex for a considerably large class of infinite matrices $M$ for all $p>1$. 
	A special attention is given to the identification of the dual space $(\ell_M^p )^*$. Building on the earlier works of Bennett and J\"agers, we extend and apply some classical factorization results to the sequence spaces $\ell_M^p$. 
	\end{abstract}
	
	
\section{Introduction}
One of the prize problems asked by the dutch mathematical society \cite{dms} was to find the dual of Ces\'aro sequence space. J\"agers \cite{Jaegers} provided an isomorphic identification, though for the larger class of spaces of sequences $(x_n)$ satisfying
\begin{equation}
\label{eq:Jaegers}
\sum_{n=1}^{\infty}(\beta_n\,\sum_{k=1}^n|x_k|)^p<\infty, \quad \text{where}\; p\geq 1\;\text{is fixed}
\end{equation}
and $(\beta_n)$ is an arbitrary sequence of positive numbers. When $\beta_n=1/n$ for all $n\in\mathbb N$, one recovers the Ces\'aro sequence space. A lot of interest is shown in the study of sequence spaces, e.g. see \cite{Boos, Wilansky}, and in particular in the Ces\'aro sequence space \cite{Yee, Lee--Ng, Lee--Ng2, Shiue1, Leibowitz, Albanese} or its counterpart Ces\'aro function space \cite{Shiue2, Maligranda1, Maligranda2, Maligranda3}. In connection with the problem of identifying the dual of the Ces\'aro sequence space, Bennett \cite{Bennett} motivated by the classical inequalities of Hilbert, Hardy, and Copson,  systematically studies the structure of the Ces\'aro sequence space and certain generalizations by employing his method of factorization of spaces. Bennett was able to give an isometric isomorphic identification for the dual of the Ces\' aro sequence space.
Further  results related to  Bennett's factorization theorems were obtained in a series of papers by \cite[Leindler]{Leindler1, Leindler2, Leindler3}. 

In view of these developments, given an infinite matrix $M=(m_{nk})$ with possibly complex-valued entries, we study the space $\ell_M^p$ of sequences $(x_n)\in\mathbb C^{\mathbb N}$ that satisfy
\begin{equation}
\label{eq:lpM}
\|x\|_{M,p}\coloneqq \Big(\sum_{n=1}^{\infty}\Big(\sum_{k=1}^{\infty}|m_{nk}|\,|x_k|\Big)^p\Big)^{1/p}<\infty.
\end{equation}
It is immediate that \eqref{eq:Jaegers} is a specialization of \eqref{eq:lpM} with $m_{nk}=\beta_n$ for every $k\leq n$ for all $n\in\mathbb N$ and $m_{nk}=0$ otherwise. These sequence spaces for $M$ with nonnegative entries were first introduced in \cite[\S17, pp.90]{Bennett}.
There are several reasons why we investigate the sequence spaces $\ell^p_M$. First they offer a natural generalization of the Ces\'aro sequence space, to which a great deal of attention has been given over the years. This generalization serves as a tool to understand better the relationship between the space $\ell^p_M$ and its underlying matrix $M$  by means of formula \eqref{eq:lpM}. The theory on infinite matrices is old, rich and it has many important applications, in particular with regard to infinite system of linear equations, e.g. see \cite[\S1.2, \S3]{Cooke}. In this respect the study of $\ell^p_M$ spaces provides a bridge between infinite matrices and the theory of sequence spaces. Secondly when studying a sequence space it is of interest to identify its Banach dual. The factorization method developed by Bennett, with minor modifications, applies directly to spaces $\ell^p_M$, thus enabling us to isometric isomorphically identify its dual for many infinite matrices $M$. Thirdly  to each infinite matrix there corresponds a matrix summability method. Summability methods are well studied and have many applications, e.g. see \cite{Boos}, therefore investigating $\ell^p_M$ spaces opens the door to possibly interesting connections with summability theory. Lastly the  Ces\'aro sequence space has seen recent applications in the theory of Dirichlet series and their multiplier algebras, e.g. see \cite{Bueno} and references therein, however in the present note we shall not deal with it.

Our work is organised as follows. In Section \ref{s:preliminaries} we give some preliminary definitions and results that are useful for the development of our work. In Section \ref{s:first-results} we present some basic results about the $\ell^p_M$ spaces. By standard methods we show that $(\ell^p_M,\|\cdot\|_{M,p})$ is a separable Banach space for all $p\geq 1$ iff $M$ has no vanishing columns (Theorem \ref{th:M-Banach-space}). In Section \ref{s:convexity} we look at the geometry of the space $\ell^p_M$ and in particular its convexity properties. We show that $(\ell^p_M,\|\cdot\|_{M,p})$  is strictly convex  for every $p>1$, if $M$ is lower triangular with nonzero diagonal terms. Moreover we prove that $(\ell^p_M,\|\cdot\|_{M,p})$ is uniformly convex for every $p>1$, whenever $M$ is invertible and $M^{-1}(\phi)\subseteq\ell^p_M$. Here $\phi$ denotes the set of all sequences in $\mathbb C^{\mathbb N}$ which are eventually vanishing, i.e. 'finite sequences' (Theorem \ref{th:bijective-M}). Then as a by-product of the well-known Milman--Pettis theorem, e.g. see any \cite{Kakutani, Milman, Pettis}, we obtain that in this case $\ell^p_M$ is reflexive.
In Section \ref{s:dual} we draw our attention to the study of the dual $(\ell^p_M)^*$. We show that for any invertible matrix $M$ the inclusion  holds $\ell^q_{(M^{-1})^T}\subseteq(\ell^p_M)^*$. If additionally $M$ is diagonal, then the reverse inclusion is also true and in particular $\ell^p_M$ is reflexive for all $p>1$ (Theorem \ref{th:invertible}). In the second part of this section we make an excursion into Bennett's factorization theorems, which with minor modifications apply to the auxiliary spaces $d_M(p), g_M(p)$ that we introduce in the same way as  Bennett does in \cite[\S3]{Bennett} (Theorem \ref{th:d-g}, Theorem \ref{th:dual}). Moreover for lower triangular matrices $M=(m_{nk})$ that have $\ell^p$-summable diagonal terms, satisfy  $|m_{n(k+1)}|\leq |m_{nk}|$ for all $k=1,2,\cdots, n-1$ for all $n\in\mathbb N$, and the terms $(m_{n1})$ of the first column enjoy a certain growth condition, we show that the factorization holds $\ell^p_M=\ell^p\cdot g_M(q)$, where $p^{-1}+q^{-1}=1$. This factorization helps us establish that the dual $(\ell^p_M)^*$ is isometric isomorphic to $d_M(q)$ (Theorem \ref{th:converse}). We illustrate most of our results with examples.

\section{Preliminaries}
\label{s:preliminaries}
\subsection{Infinite matrices}
Let $M=(m_{nk})_{n,k\in\mathbb N_0}$ be an infinite matrix with possibly complex-valued entries. Given an element $x\in\mathbb C^{\mathbb N}$ define
\begin{equation}
\label{eq:M-x}
(Mx)_n\coloneqq\sum_{k=1}^{\infty}m_{nk}\,x_k\quad (n=1,2,3,\cdots).
\end{equation}
The inverse of a matrix $M$, if it exists, is a matrix $M^{-1}$ such that $MM^{-1}=M^{-1}M=I$. 
By $\phi$ we denote the set of all finite sequences in $\mathbb C^{\mathbb N}$, i.e. $x=(x_n)\in\phi$ if $x_n=0$ for all $n\geq n_0$ for some $n_0\in\mathbb N$. Evidently if $x\in\phi$, then  $(Mx)_n$ exists for every $n\in\mathbb N$.
\subsubsection{Some examples}
\begin{itemize}
\item the identity matrix $I=(m_{nk})$, where $m_{nk}=1$ if $n=k$ and $m_{nk}=0$ otherwise.
\item the generalized Ces\'aro matrix $M=C_{\alpha},\,\alpha\in\mathbb R\setminus\{-\mathbb N\}$, e.g. see \cite[\S 3, pp.104]{Boos}, with general term 
\begin{equation}
\label{eq:Cesaro}
m_{nk}=\left\{
\begin{array}{ll}
      \frac{\displaystyle{n-k+\alpha-1\choose n-k}}{\displaystyle{n+\alpha-1\choose n-1}} & k\leq n \\[1em]
     0 & \text{else}.
\end{array}
\right.
\end{equation}
\item N\"orlund matrix $N=(m_{nk})$ where $m_{nk}=p_{n-k}/P_n$ for $k\leq n$ and $m_{nk}=0$ else and $(p_k)$ is a sequence of positive numbers and $P_n=\sum_{k=0}^{n}p_k$ . This includes in particular Zweier method (of order $\alpha\neq 0$), e.g. see \cite[\S3, pp.127]{Boos}.
\item Riesz matrix $R=(m_{nk})$ where $m_{nk}=p_{k}/P_n$ for $k\leq n$ and $m_{nk}=0$ else and $(p_k)$ is a sequence of positive numbers and $P_n=\sum_{k=0}^{n}p_k$.
\item Hausdorff matrix $H=(m_{nk})$, e.g. see \cite[\S3, pp. 137]{Boos}, with general term $m_{nk}=\Delta\text{diag}(p_n)\Delta$ where $\text{diag}(p_n)$ is a diagonal matrix determined by $(p_n)\in\mathbb C^{\mathbb N}$ and $\Delta$ is the lower triangular matrix with general term $\Delta_{nk}=(-1)^k{n\choose k}$.
\item Hilbert matrix where $\mathcal H=(m_{nk})$ is given by $m_{nk}=1/(n+k-1)$, e.g. see \cite[pp.53]{Bennett}.
\end{itemize}
\subsection{Strict and uniform convexity}
A normed space $(X,\|\cdot\|)$ is strictly convex if for any $x,y\in X$ with $\|x\|=\|y\|=1$ we have $\frac{1}{2}\|x+y\|<1$. A stronger notion of convexity is that of uniform convexity. A normed space $(X,\|\cdot\|)$ is uniformly convex  if and only if for every $0<\varepsilon\leq 2$ there is $\delta(\varepsilon)>0$ such that $\frac{1}{2}\|x+y\|\leq 1-\delta(\varepsilon)$ whenever $\|x\|=\|y\|=1$ and $\|x-y\|\geq \varepsilon$. For $0\leq \varepsilon\leq 2$ the modulus of convexity of $X$ is defined as
\begin{equation}
\label{eq:modulus-convex}
\delta_X(\varepsilon):=\inf\{1-\frac{1}{2}\,\|x+y\|\,:\,\|x\|=\|y\|=1,\;\|x-y\|\geq \varepsilon\}.
\end{equation}
Evidently $X$ is uniformly convex if and only if $\delta_X(\varepsilon)>0$ for all $\varepsilon\in(0,2]$. We consider another quantity that measures uniform convexity (Gurari\'i's modulus of convexity). For $0\leq \varepsilon\leq 1$ let 
\begin{equation}
\label{eq:Gurarii}
\beta_X(\varepsilon):=\inf\{\sup_{0\leq\alpha\leq 1}(1-\|\alpha x+(1-\alpha)y\|)\,:\,\|x\|=\|y\|=1,\,\|x-y\|\geq \varepsilon\}.
\end{equation}
It was shown by Gurari\'i, e.g. see \cite[Theorem 1]{Gurarii}, that
\begin{equation}
\label{eq:Gurarii-unif-convex}
\delta_X(\varepsilon)\leq \beta_X(\varepsilon)\leq 2\,\delta_X(2\varepsilon).
\end{equation}
It follows that $X$ is uniformly convex if and only if $\beta_X(\varepsilon)>0$ for every $\varepsilon\in(0,1]$. 

\subsection{Banach spaces}
A complete normed linear space $(X,\|\cdot\|)$ is called a Banach space. By $(X^*,\|\cdot\|_*)$ we denote its dual, the space of all bounded linear functionals on $X$. The dual $X^*$ is itself a Banach space with its norm defined by 
\begin{equation}
\label{eq:dual-norm}
\|f\|_*:=\sup_{x\in X, x\neq 0}\frac{|f(x)|}{\|x\|},\quad \text{for all}\;f\in X^*.
\end{equation}
Denote by $(X^{**},\|\cdot\|_{**})$ the bidual of $X$, that is the dual of the dual $X^*$. For each $x\in X$ let $J(x):X^*\to\mathbb C$ be the evaluation scalar map generated by $x$ via the rule $J(x)(f)=f(x)$. It is known that $J(x)$ is an injective continuous linear functional on $X^*$, i.e. $J(x)\in X^{**}$ and that it preserves norms, i.e. for every $x\in X$ it holds that $\|J(x)\|_{**}=\|x\|$. A Banach space $(X,\|\cdot\|)$ is reflexive if $J(x)$ is surjective for every $x\in X$. 
\begin{mainthm}\cite{Milman}
\label{th:Milman}
Uniformly convex Banach spaces are reflexive.
\end{mainthm}
A sequence of elements $(v_n)\subseteq X$ is a Schauder basis for $X$ if for every $x\in X$ there is a unique sequence $(a_n)\in\mathbb C^{\mathbb N}$ such that 
\begin{equation}
\label{eq:Schauder}
x=\sum_{n=1}^{\infty}a_n\,v_n.
\end{equation}
Equation \eqref{eq:Schauder} means that the series converges in the chosen order of the elements. If the order plays no role for any $x\in X$ then the basis $(v_n)$ is unconditional.  Note that the definition requires that the sequence $(v_n)$ be complete in $X$.

\begin{lemma}\cite[Theorem 3.1.4]{Christensen}
\label{l:Schauder}
A complete  sequence of elements $(v_n)\subseteq X$ is a Schauder basis for $X$ if and only if there exists a constant $K>0$ such that for all $m,n\in\mathbb N$ with $m\leq n$ the inequality holds
\begin{equation}
\label{eq:Schauder-inequality}
\|\sum_{k=1}^ma_k\,v_k\|\leq K\,\|\sum_{k=1}^na_k\,v_k\|,\quad\text{for all sequences}\;(a_k)\in\mathbb C^{\mathbb N}.
\end{equation}
\end{lemma}



\section{First results for $\ell^p_M$ spaces}
\label{s:first-results}
\subsection{Sequence spaces $\ell^p_M$}
For $x\in\mathbb C^{\mathbb N}$ and $p\geq 1$ we let 
\begin{equation}
\label{eq:norm}
\|x\|_{M,p}\coloneqq \Big(\sum_{n=1}^{\infty}\Big(\sum_{k=1}^{\infty}|m_{nk}|\,|x_k|\Big)^p\Big)^{1/p}.
\end{equation}
Given a matrix $M$ we define the sequence space
\begin{equation}
\label{eq:seq-space}
\ell^p_M\coloneqq\{x\in\mathbb C^{\mathbb N}\,:\, \|x\|_{M,p}<\infty\}.
\end{equation}

\begin{proposition}
\label{p:non-van-col}
$\|\cdot\|_{M,p}$ is a norm in $\ell^p_M$ if and only if $M$ has no vanishing columns. 
\end{proposition}
\begin{proof}
Suppose there is $k\in\mathbb N$ such that $m_{nk}=0$ for all $n\in\mathbb N$. Consider $x=(x_j)_{j\in\mathbb N}\in\mathbb C^{\mathbb N}$ with $x_j\neq 0$ when $j=k$ and $x_j=0$ for all $j\neq k$. Then $(Mx)_n=\sum_{j=1}^{\infty}m_{nj}\,x_j=m_{nk}\,x_k=0$ for all $n\in\mathbb N$ implies $\|x\|_{M,p}=(\sum_{n=1}^{\infty}|m_{nk}|^p\,|x_k|^p) ^{1/p}=0$ while $x\neq 0$. Consequently $\|\cdot\|_{M,p}$ is not positive definite and thus not a norm. 

Now let $M$ have no vanishing columns. It is clear that $\|\cdot\|_{M,p}$ is non-negative and for any $\alpha\in\mathbb C$ we have $\|\alpha\,x\|_{M,p}=|\alpha|\,\|x\|_{M,p}$. Moreover by Minkowski's inequality it follows that $\|x+y\|_{M,p}\leq \|x\|_{M,p}+\|y\|_{M,p}$ for any $x,y\in\ell^p_M$. Let $\|x\|_{M,p}=0$, then $|m_{nk}|\,|x_k|=0$ for all $k,n\in\mathbb N$. But for every $k\in\mathbb N$ there is $n(k)\in\mathbb N$ such that $m_{n(k)k}\neq 0$ implying $x_k=0$ for every $k\in\mathbb N$, i.e. $x=0$. 
\end{proof}

In view of Proposition \ref{p:non-van-col} we restrict to matrices $M$ with non-vanishing columns.

\begin{thm}
\label{th:M-Banach-space}
$(\ell^p_M,\|\cdot\|_{M,p})$ is a separable Banach space for all $p\geq 1$. 
\end{thm}
\begin{proof}
By Proposition \ref{p:non-van-col} we have that $(\ell^p_M,\|\cdot\|_{M,p})$ is a normed space. We show completeness. Let $(x_n)$ be a Cauchy sequence in $\ell^p_M$. 
Then for every $\varepsilon>0$ there is $N(\varepsilon)\in\mathbb N$ such that $\|x_n-x_m\|_{M,p}<\varepsilon$ for all $m,n\geq N(\varepsilon)$. This in particular implies 
\begin{equation}
\label{eq:Cauchy}
|m_{kj}|\,|x_{nj}-x_{mj}|<\varepsilon,\quad\text{for all}\;m,n\geq N(\varepsilon),\;j,k\in\mathbb N.
\end{equation}
Because $M$ has non-vanishing columns then for every $j\in\mathbb N$ there is $k(j)\in\mathbb N$ such that $|m_{k(j)j}|>0$, then 
\begin{equation}
\label{eq:Cauchy2}
|x_{nj}-x_{mj}|<\frac{\varepsilon}{|m_{k(j)j}|},\quad\text{for all}\;m,n\geq N(\varepsilon),\;j\in\mathbb N.
\end{equation}
Hence $(x_{nj})$ is a Cauchy sequence in $\mathbb C$ for every $j\in\mathbb N$ consequently $x_{nj}\to x_{j}$ for a certain $x_{j}\in\mathbb C$ as $n\to\infty$ for every $j\in\mathbb N$.
For fixed $J,K\in\mathbb N$ we have
\begin{equation}
\label{eq:finite-J}
\Big(\sum_{k=1}^K\Big(\sum_{j=1}^J|m_{kj}|\,|x_{nj}-x_{j}|\Big)^p\Big)^{1/p}<\varepsilon,\quad \text{for all}\;n\geq N(\varepsilon).
\end{equation}
By Minkowski's inequality we get 
\begin{align*}
\Big(\sum_{k=1}^K\Big(\sum_{j=1}^J|m_{kj}|\,|x_{j}|\Big)^p\Big)^{1/p}
<\varepsilon +\Big(\sum_{k=1}^K\Big(\sum_{j=1}^J|m_{kj}|\,|x_{nj}|\Big)^p\Big)^{1/p}
\end{align*}
for all $n\geq N(\varepsilon)$. Letting $J, K\to\infty$ yields
$$\|x\|_{M,p}=\Big(\sum_{k=1}^{\infty}\Big(\sum_{j=1}^{\infty}|m_{kj}|\,|x_j|\Big)^p\Big)^{1/p}\leq \varepsilon+\Big(\sum_{k=1}^{\infty}\Big(\sum_{j=1}^{\infty}|m_{kj}|\,|x_{nj}|\Big)^p\Big)^{1/p}=\varepsilon+\|x_n\|_{M,p}$$
for all $n\geq N(\varepsilon)$ therefore $x\in \mathscr \ell^p_M$ and by \eqref{eq:finite-J} letting $J,K\to\infty$ gives $x_n\to x$.  
Let $x^N=(x_k^N)_{k\leq N}$ be the truncation of $x\in\mathbb C^{\mathbb N}$ where $N\in\mathbb N$. Then for every $\varepsilon>0$ there is $N(\varepsilon)\in\mathbb N$ such that $\|x-x^N\|_{M,p}<\varepsilon/2$ and there is $x^N_{\mathbb Q}\in(\mathbb Q+i\,\mathbb Q)^N$ such that $\|x^N-x^N_{\mathbb Q}\|_{M,p}<\varepsilon/2$ implying that $\|x-x^N_{\mathbb Q}\|_{M,p}<\varepsilon$, thus $\ell^p_M$ is separable. 
\end{proof}

\begin{proposition}
\label{p:basic}
It holds $\ell^p_M=\{0\}$ if and only if $\sum_{n=1}^{\infty}|m_{nk}|^p=\infty$ for all $k\in\mathbb N$. Moreover $\ell^p_M\subset\ell^q_M$ whenever $p< q$. 
\end{proposition}
\begin{proof}
Let $\ell^p_M=\{0\}$ then $x\in\mathbb C^{\mathbb N}$ with $x_n=1$ for $n=k$ and $x_n=0$ else is not an element of  $\ell^p_M$ implying $\sum_{n=1}^{\infty}|m_{nk}|^p=\sum_{n=1}^{\infty}(\sum_{j=1}^{\infty}|m_{nj}|\,|x_j|)^p=\|x\|_M^p=\infty$ for every $k\in\mathbb N$. Conversely if $\sum_{n=1}^{\infty}|m_{nk}|^p=\infty$ for all $k\in\mathbb N$ and $x\neq 0$ then for some $k_0\in\mathbb N$ we would have $x_{k_0}\neq 0$ and $\|x\|^p_M=\sum_{n=1}^{\infty}(\sum_{j=1}^{\infty}|m_{nj}|\,|x_j|)^p=\sum_{n=1}^{\infty}|m_{nk_0}|^p=\infty$ implying $x\notin \ell^p_M$. Hence $\ell^p_M=\{0\}$. Next we prove the second assertion. Let $p<q$  and $x\in\ell^p_M$ then $\|x\|_{M,p}<\infty$ implies in particular that $\sum_{k=1}^{\infty}|m_{nk}|\,|x_k|<1$ for all sufficiently large $n\in\mathbb N$ and so $(\sum_{k=1}^{\infty}|m_{nk}|\,|x_k|)^q<(\sum_{k=1}^{\infty}|m_{nk}|\,|x_k|)^p$ whenever $q>p$ for all sufficiently large $n\in\mathbb N$. Then evidently $\|x\|_{M,q}<\infty$, i.e. $x\in\ell^q_M$.  
\end{proof}

\begin{proposition}
\label{p:Schauder-basis} 
The sequence $v_n=(0,\cdots, 0, \underbrace{1}_{n-th\;term} , 0, \cdots)$  for $n\in\mathbb N$ is an unconditional Schauder basis for every $(\ell^p_M,\|\cdot\|_{M,p})$ for all $p\geq 1$. 
\end{proposition}
\begin{proof}
It follows as an application of Lemma \ref{l:Schauder}. First we prove that $(v_n)$ is complete. Let $V=\{x\in\ell^p_M\,:\,x=\sum_{k=1}^{n}a_k\,v_k\;\text{for some not all zero}\;a_k\in\mathbb C,\;n\in\mathbb N\}$, then we show that $\ell^p_M=\cl V$. Clearly $V\subseteq\ell^p_M$. 
Now let $x\in\ell^p_M$ and denote by $x^N$ the truncation of $x$, i.e. $x^N_k=x_k$ for every $k\leq N$ and $x^N_k=0$ for all $k>N$, then $x^N\in V$ for every $N\in\mathbb N$. Moreover 
\begin{align*}
\|x-x^N\|_{M,p}=(\sum_{n=1}^{\infty}(\sum_{k\geq N+1}|m_{nk}|\,|x_k|)^p)^{1/p}\to 0\;\text{as}\;N\to\infty
\end{align*}
since $\sum_{k\geq N+1}|m_{nk}|\,|x_k|\to 0$ as $N\to\infty$ for all $n\in\mathbb N$ by Cauchy's criterion. This proves $\ell^p_M\subseteq\cl V$. 
Next let $\sigma$ be a permutation of $\mathbb N$, take $(a_k)\in\mathbb C^{\mathbb N}$ and $i,j\in\mathbb N$ with $i\leq j$, then 
\begin{align*}
\|\sum_{k=1}^ia_{k}\,v_{\sigma(k)}\|_{M,p}&=(\sum_{n=1}^{\infty}(\sum_{k=1}^i|m_{n\sigma(k)}|\,|a_k|)^p)^{1/p}\\&\leq (\sum_{n=1}^{\infty}(\sum_{k=1}^j|m_{n\sigma(k)}|\,|a_k|)^p)^{1/p}=\|\sum_{k=1}^ja_k\,v_{\sigma(k)}\|_{M,p}
\end{align*}
implies that inequality \eqref{eq:Schauder-inequality} holds with $K=1$ for every permutation $\sigma$ of $\mathbb N$.
\end{proof}


\subsection{Convexity of $\ell^p_M$}
\label{s:convexity}
\begin{thm}
\label{th:bijective-M}
The followings are true:
\begin{enumerate}[(a)]
\item If $M$ is lower triangular with nonzero diagonal terms and $p>1$, then $(\ell^p_M,\|\cdot\|_{M,p})$ is strictly convex.
\item If $M=(m_{nk})$ is invertible with inverse $M^{-1}$ such that $M^{-1}(\phi)\subset\ell^p_M$, then $(\ell^p_M,\|\cdot\|_{M,p})$ is uniformly convex for every $p>1$.
\end{enumerate}
\end{thm}
\begin{proof}
\begin{enumerate}[(a)]
\item Let $M$ be lower triangular matrix with  with $m_{nn}\neq 0$ and let $p>1$.
Take $x,y\in\ell^p_M$ with $\|x\|_{M,p}=\|y\|_{M,p}=1$ and $\|x+y\|_{M,p}=2$. 
By Minkowski's inequality we have $2=\|x+y\|_{M,p}\leq \|x\|_{M,p}+\|y\|_{M,p}=2$. By the strict convexity of $\|\cdot\|_p$ when $p>1$ we obtain that 
$$\sum_{k=1}^{n}|m_{nk}|\,|x_k|=\sum_{k=1}^{n}|m_{nk}|\,|y_k|$$
for all $n\in\mathbb N$, implying $|x_n|=|y_n|$ for all $n\in\mathbb N$. We aim to show that $x_n=y_n$ for all $n\in\mathbb N$. Define the sets $$D=\{n\in\mathbb N\,:\,x_n\neq y_n,\,|x_n|=|y_n|\}\; \text{and}\; D_1=\{n\in\mathbb N\,:\,x_n=-y_n,\,x_n\neq 0, |x_n|=|y_n|\}.$$ Evidently $D_1\subseteq D$ and 
\begin{align*}
\sum_{k=1}^{\infty}|m_{nk}|\,\Big|\frac{x_k+y_k}{2}\Big|&=\sum_{k\in\mathbb N\setminus D}|m_{nk}|\,|x_k|+\sum_{k\in D\setminus D_1}|m_{nk}|\,\Big|\frac{x_k+y_k}{2}\Big|.
\end{align*}
Suppose that $D\neq \emptyset$. We distinguish two cases. First if $D_1\neq \emptyset$ then
\begin{align*}
\sum_{k\in\mathbb N\setminus D}|m_{nk}|\,|x_k|+\sum_{k\in D\setminus D_1}|m_{nk}|\,\Big|\frac{x_k+y_k}{2}\Big|\leq \sum_{k\in\mathbb N\setminus D_1}|m_{nk}|\,|x_k|<\sum_{k\in \mathbb N}|m_{nk}|\,|x_k|.
\end{align*}
If $D_1=\emptyset$, then $D\setminus D_1\neq\emptyset$ implying that there is $n\in\mathbb N$ such that $|x_n|=|y_n|$ but $x_n\neq \pm y_n$. If $\theta$ is the smaller angle between the vectors $\overrightarrow{0x_n}$ and $\overrightarrow{0y_n}$ in $\mathbb C$, then $\theta\in(0,\pi)$ and in particular $|x_n+y_n|<|x_n|+|y_n|$. Consequently 
\begin{align*}
\sum_{k\in\mathbb N\setminus D}|m_{nk}|\,|x_k|&+\sum_{k\in D\setminus D_1}|m_{nk}|\,\Big|\frac{x_k+y_k}{2}\Big|\\&<\sum_{k\in\mathbb N\setminus D}|m_{nk}|\,|x_k|+\sum_{k\in D\setminus D_1}|m_{nk}|\,\frac{|x_k|+|y_k|}{2}=\sum_{n\in\mathbb N}|m_{nk}|\,|x_k|.
\end{align*}
In either case we have
$$\sum_{k=1}^{\infty}|m_{nk}|\,\Big|\frac{x_k+y_k}{2}\Big|<\sum_{k=1}^{\infty}|m_{nk}|\,|x_k|$$
implying
\begin{align*}
1=\Big\|\frac{x+y}{2}\Big\|_{M,p}&=\Big(\sum_{n=1}^{\infty}\Big(\sum_{k=1}^{n}|m_{nk}|\,\Big|\frac{x_k+y_k}{2}\Big|\Big)^p\Big)^{1/p}\\&=\Big(\sum_{n=1}^{\infty}\Big(\sum_{k=1}^{\infty}|m_{nk}|\,\Big|\frac{x_k+y_k}{2}\Big|\Big)^p\Big)^{1/p}\\&
<\Big(\sum_{n=1}^{\infty}\Big(\sum_{k=1}^{\infty}|m_{nk}|\,|x_k|\Big)^p\Big)^{1/p}=\Big(\sum_{n=1}^{\infty}\Big(\sum_{k=1}^{n}|m_{nk}|\,|x_k|\Big)^p\Big)^{1/p}=\|x\|_{M,p}=1
\end{align*}
which is impossible. Therefore $D=\emptyset$, i.e. $x=y$. 

\item 
If $p=1$ by taking $M=I$ we know that $\ell^1$ is not uniformly convex, thus let $p>1$. Take $u=((1-(\frac{\varepsilon}{2})^p)^{\frac{1}{p}}, \frac{\varepsilon}{2}, 0, 0, \cdots)$ and $v=((1-(\frac{\varepsilon}{2})^p)^{\frac{1}{p}}, -\frac{\varepsilon}{2}, 0, 0, \cdots)$. Note that $\phi\subset \omega_{M^{-1}}$, therefore there are $x_0,y_0\in\mathbb C^{\mathbb N}$ such that $x_0=M^{-1}u$ and $y_0=M^{-1}v$. Moreover the assumption $M^{-1}(\phi)\subset\ell^p_M$ implies $x_0, y_0\in\ell^p_M$.
By the relation $\|z\|_{M,p}\geq \|Mz\|_p$ for every $z\in\mathbb C^{\mathbb N}$ we obtain 
\begin{align*}
&\|x_0\|_{M,p}\geq \|Mx_0\|_p=\|M(M^{-1}u)\|_p=\|u\|_p=1\\&
\|y_0\|_{M,p}\geq \|My_0\|_p=\|M(M^{-1}v)\|_p=\|v\|_p=1.
\end{align*}
Moreover by similar calculations $\|x_0-y_0\|_{M,p}\geq \|M(x_0-y_0)\|_p=\|u-v\|_p=\varepsilon$. 
Define the auxiliary quantities
\begin{align*}
&\widehat\beta_{\ell^p_M}(\varepsilon):=\inf\{\sup_{0\leq\alpha\leq 1}(1-\|\alpha x+(1-\alpha)y\|_{M,p})\,:\,\|x\|_{M,p}\geq 1,\,\|y\|_{M,p}\geq 1,\,\|x-y\|_{M,p}\geq \varepsilon\}\\&
\widehat{\delta}_{\ell^p_M}(\varepsilon):=\inf\{1-\frac{1}{2}\|x+y\|_{M,p}\,:\,\|x\|_{M,p}\geq 1,\,\|y\|_{M,p}\geq 1,\,\|x-y\|_{M,p}\geq \varepsilon\}.
\end{align*}
Evidently $\widehat{\beta}_{\ell^p_M}(\varepsilon)\geq \widehat{\delta}_{\ell^p_M}(\varepsilon)$ for every $\varepsilon\in[0,1]$. From \eqref{eq:Gurarii} we have that $\beta_{\ell^p_M}(\varepsilon)\geq \widehat{\beta}_{\ell^p_M}(\varepsilon)$ for every $\varepsilon\in[0,1]$.  Therefore in view of \eqref{eq:Gurarii-unif-convex} it suffices to prove that $\widehat{\delta}_{\ell^p_M}(\varepsilon)>0$ for every $\varepsilon\in(0,1]$. For our choice of $x_0$ and $y_0$ it holds that 
\begin{align*}
\sup_{0\leq\alpha\leq 1}(1-\|\alpha x_0+(1-\alpha)y_0\|_{M,p})&=1-\inf_{0\leq\alpha\leq 1}\|\alpha x_0+(1-\alpha)y_0\|_{M,p}\\&\leq 1-\inf_{0\leq \alpha\leq 1}\|\alpha u+(1-\alpha)v\|_p=1-\Big(1-\Big(\frac{\varepsilon}{2}\Big)^p\Big)^{\frac{1}{p}}
\end{align*}
consequently $$\widehat{\beta}_{\ell^p_M}(\varepsilon)\leq 1-\Big(1-\Big(\frac{\varepsilon}{2}\Big)^p\Big)^{\frac{1}{p}}$$
implying
\begin{align*}
\Big(1-\Big(\frac{\varepsilon}{2}\Big)^p\Big)^{\frac{1}{p}}\leq 1-\widehat{\beta}_{\ell^p_M}(\varepsilon)\leq 1-\widehat{\delta}_{\ell^p_M}(\varepsilon),\quad \varepsilon\in[0,1].
\end{align*}
The quantity on the left is strictly positive and $<1$ for every $\varepsilon\in(0,1]$, thus $\widehat{\delta}_{\ell^p_M}(\varepsilon)>0$ for every $\varepsilon\in(0,1]$. This completes the proof.
\end{enumerate}
\end{proof}

As a direct consequence of Milman--Pettis Theorem  we obtain that:

\begin{corollary}
\label{c:reflexive} 
Let $M$ satisfy conditions of Theorem \ref{th:bijective-M}(b), then $(\ell^p_M,\|\cdot\|_{M,p})$ is reflexive for every $p>1$.
\end{corollary}

\begin{example}
Consider the standard Ces\'aro matrix $M=C$ where $m_{nk}=1/n$ for $k\leq n$ and $m_{nk}=0$ else, for every $n\in\mathbb N$. Then $C$ is invertible with inverse
\begin{align*}
C^{-1}=\begin{pmatrix}
\phantom{-}1 & \phantom{-}0 & \phantom{-}0 & 0 & \cdots\\
-1 & \phantom{-}2 & \phantom{-}0 & 0 &\cdots\\
\phantom{-}0 & -2 & \phantom{-}3 & 0 & \cdots\\
\phantom{-}0 & \phantom{-}0 & -3 & 4 & \cdots\\
\cdots &\cdots & \cdots & \cdots &\cdots
\end{pmatrix}.
\end{align*}
Let $u=(u_n)\in\phi$, then $u_n=0$ for all $n\geq m$ for some $m\in\mathbb N$. Consequently 
\begin{align*}
C^{-1}u=\begin{pmatrix}
u_1\\ -u_1+2u_2\\ -2u_2+3u_3\\\cdots\\-(m-1)u_{m-1}+mu_m\\ -mu_m\\0\\0\\\cdots
\end{pmatrix}=\begin{pmatrix}
x_1\\ x_2\\ x_3\\\cdots\\x_{m}\\x_{m+1}\\0\\0\\\cdots
\end{pmatrix}=x\in\phi.
\end{align*} 
This implies $C^{-1}(\phi)\subset\phi$.
By a theorem of Bennett when $p>1$, e.g. see \cite[Theorem 1.5]{Bennett}, it holds that $\phi\subset\ell^p_C$ (commonly written $\text{ces}(p)$) and in particular $C^{-1}(\phi)\subset\ell^p_C$. Conditions of Theorem \ref{th:bijective-M} are fulfilled implying that $\ell^p_C$ is uniformly convex for all $p>1$.
\end{example}


\section{The dual of $\ell^p_M$} 
\label{s:dual}

\subsection{The case of invertible matrices}

\begin{thm}
\label{th:invertible}
Let $M=(m_{nk})$ be an invertible matrix, then $\ell^q_{(M^{-1})^T}\subseteq (\ell^p_M)^*$. If $M$ is also diagonal then the reverse inclusion holds, in particular   $\ell^p_M$ is reflexive for $p>1$.
\end{thm}
\begin{proof}
Let $M$ have an inverse $M^{-1}=(m^*_{nk})$. We consider the case $p>1$. When $p=1$ similar steps apply. Take $y\in\ell^q_{(M^{-1})^T}$ and $x\in\ell^p_M$. For $N\in\mathbb N$ let $y^N, x^N$ be the truncation of $y$ and $x$ respectively up to the $N$-th term i.e. $x^N_n=x_n,\,y^N_n=y_n$ for $n\leq N$ and $x^N_n=y^N_n=0$ for $n>N$. Note that $x^N\in\ell^p_M,\;y^N\in\ell^q_{(M^{-1})^T}$ and
\begin{align*}
\sum_{n=1}^{N}y_n\,x_n&=\sum_{n=1}^{N}y_n\,(M^{-1}(Mx))_n=\sum_{n=1}^{N}y_n\sum_{k=1}^{\infty}m^*_{nk}\,(Mx^N)_k\\&
=\sum_{k=1}^{\infty}(Mx^N)_k\,\sum_{n=1}^Nm^*_{nk}\,y^N_n=\sum_{k=1}^{\infty}(Mx^N)_k\,((M^{-1})^Ty^N)_k
\\&\leq \Big(\sum_{k=1}^{\infty}|((M^{-1})^Ty^N)_k|^q\Big)^{1/q}\,\Big(\sum_{k=1}^{\infty}|(Mx^N)_k|^p\Big)^{1/p}
= \|y^N\|_{(M^{-1})^T,q}\,\|x^N\|_{M,p}.
\end{align*}
Letting $N\to\infty$ we then obtain 
$$\sum_{n=1}^{\infty}y_n\,x_n\leq \|y\|_{(M^{-1})^T,q}\,\|x\|_{M,p}$$
implying $\|y\|_*\leq \|y\|_{(M^{-1})^T,q}<\infty$ therefore $y\in(\ell^p_M)^*$.

Now suppose that additionally $M=(m_{nk})$ is a diagonal matrix. Because $M$ is invertible then $m_{nn}\neq 0$ for all $n\in\mathbb N$. In particular the inverse $M^{-1}$ is given by the diagonal matrix $M^{-1}$ with entries $1/m_{nn}$ for all $n\in\mathbb N$. Note that $(M^{-1})^T=M^{-1}$. Now let $y\in(\ell^p_M)^*$ and take $x\in\mathbb C^{\mathbb N}$ defined by 
\begin{equation*}
x_n\coloneq \left\{
\begin{array}{ll}
      \displaystyle \frac{|y_n|^{q-2}}{|m_{nn}|^q}\,\overline{y}_n & n\leq N \\[1em]
     0 & n>N.
\end{array} 
\right.
\end{equation*}
Here $\overline{y}_n$ is the complex conjugate of $y_n$ for $n\in\mathbb N$.
Then we get that $x\in\ell^p_M$, since all but finitely many terms are nonzero, with norm given by
\begin{align*}
\|x\|_{M,p}=(\sum_{n=1}^{\infty}(|m_{nn}|\,|x_n|)^p)^{1/p}=\Big(\sum_{n=1}^{N}\frac{|y_n|^q}{|m_{nn}|^q}\Big)^{1/p}.
\end{align*}
On the other hand we have
\begin{align*}
\sum_{n=1}^{\infty}y_n\,x_n=\sum_{n=1}^{N}\frac{|y_n|^q}{|m_{nn}|^q}\leq \|y\|_*\,\|x\|_{M,p}
\end{align*}
implying 
\begin{align*}
\|y\|_*\geq \Big(\sum_{n=1}^{N}\frac{|y_n|^q}{|m_{nn}|^q}\Big)^{1/q},\quad N\in\mathbb N.
\end{align*}
Therefore $\|y\|_{M^{-1},q}\leq\|y\|_*<\infty$, consequently $y\in\ell^q_{M^{-1}}$. 
Following the same line of arguments we can show that $(\ell^q_{M^{-1}})^*=\ell^p_M$ since $(M^{-1})^{-1}=M$. In particular $\ell^p_M=(\ell^p_M)^{**}$, hence $\ell^p_M$ is reflexive. This completes the proof.
\end{proof}

\begin{remark}
Reflexivity follows immediately as a consequence of \cite[Proposition 17.18]{Bennett} because any lower-triangular matrix $M$ fulfills the finite-row condition.
\end{remark}

\begin{remark}
\label{r:counter-example}
The reverse inclusion is not in general true for a non-diagonal matrix. Take $M=(m_{nk})$ given by $m_{n(n+1)}=m_{nn}=1$ and $m_{nk}=0$ else, then $M^{-1}=(m^*_{nk})$ with $m^*_{nk}=(-1)^{n+k}$ for $k\geq n$ and $m_{nk}=0$ for $k<n$ for every $n\in\mathbb N$. Note that $M$ has $\ell^p$-summable columns, but $(M^{-1})^{T}$ has no $\ell^q$-summable columns. In view of Proposition \ref{p:basic} then  $\ell^p_M\neq \{0\}$ and $\ell^q_{(M^{-1})^T}=\{0\}$. But by Hahn--Banach Theorem, e.g. see \cite[Theorem 1, pp. 106]{Yosida}, the dual space $(\ell^p_M)^*\neq\{0\}$ as $\ell^p_{M}\neq\{0\}$.
\end{remark}

\begin{example}
Let $C=(m_{nk})$ be the Ces\'aro method. We know that $C$ is injective and thus it has an inverse $C^{-1}=(m^*_{nk})$ given by $m^*_{n(n-1)}=-(n-1),\,m^*_{nn}=n$ for all $n\in\mathbb N$ and $m_{nk}=0$ otherwise. Note that $(C^{-1})^T$ has $\ell^q$-summable columns and so in view of Proposition \ref{p:basic} it follows that $\ell^q_{(C^{-1})^T}\neq\{0\}$. On the other side we have that $(\ell^p_C)^*$ is isomorphic to the space $d(q)$ and isometric to it when $\ell^p_C$ is equipped with a certain equivalent norm, e.g. see \cite[Corollary 12.17]{Bennett}. So we have $\{0\}\subsetneq \ell^q_{(C^{-1})^T}\subseteq (\ell^p_C)^*\simeq d(q)$. 
\end{example}


\subsection{Two essential lemmas} The next lemma plays a key role in several places. It appears as \cite[Lemma 3.6]{Bennett} and we include its proof for completeness.
\begin{lemma}[Summation by parts]
\label{l:partial-sums}
Let $(u_n),(v_n),(w_n)$ be nonnegative sequences such that $(w_n)$ is nonincreasing then 
\begin{equation}
\label{eq:first-ineq}
\sum_{k=1}^nu_k\leq \sum_{k=1}^nv_k\quad (n=1,2,3,\cdots)
\end{equation}
implies
\begin{equation}
\label{eq:second-ineq}
\sum_{k=1}^{n}u_k\,w_k\leq\sum_{k=1}^nv_k\,w_k\quad (n=1,2,3,\cdots).
\end{equation}
\end{lemma}
\begin{proof}
Let $U_n=\sum_{k=1}^nu_k,\,V_n=\sum_{k=1}^nv_k$ then 
\begin{align*}
\sum_{k=1}^nu_k\,w_k=U_n\,w_n-\sum_{k=1}^{n-1}U_k(w_{k+1}-w_k).
\end{align*}
On the other hand $U_n\leq V_n$ and $w_{n+1}\leq w_n$ for all $n\in\mathbb N$ implies 
\begin{align*}
\sum_{k=1}^{n-1}U_k\,(w_{k+1}-w_k)\geq \sum_{k=1}^{n-1}V_k\,(w_{k+1}-w_k)
\end{align*}
therefore 
\begin{align*}
\sum_{k=1}^nu_k\,w_k\leq V_n\,w_n-\sum_{k=1}^{n-1}V_k(w_{k+1}-w_k).
\end{align*}
Moreover we have 
\begin{align*}
\sum_{k=1}^{n-1}V_k(w_{k+1}-w_k)=V_{n-1}\,w_n-V_1\,w_1-\sum_{k=2}^{n-1}w_k\,(V_k-V_{k-1}).
\end{align*}
Realizing $V_k-V_{k-1}=v_k$ and $V_1=v_1$ we obtain
\begin{align*}
\sum_{k=1}^nu_k\,w_k&\leq (V_n-V_{n-1})\,w_n+V_1\,w_1+\sum_{k=2}^{n-1}w_k\,(V_k-V_{k-1})\\&
\leq v_n\,w_n+v_1\,w_1+\sum_{k=2}^{n-1}w_k\,v_k=\sum_{k=1}^nv_k\,w_k.
\end{align*}
\end{proof}

\begin{lemma}[Bennett's partition lemma]
\label{l:Bennett}
Let $p\geq 1$ and $x\in\ell^p$. Let $(i_n)_{n\in\mathbb N_0}$ be a sequence defined recursively as follows, $i_0=0$ and with general term
\begin{equation}
\label{eq:ratios}
i_{n+1}=\sup\Big\{t\in\mathbb N,\,t>i_n\,:\,\frac{|x_{t}|^p+\cdots+|x_{i_n+1}|^p}{a_{t}(p)+\cdots+a_{i_n+1}(p)}=\sup_{\widetilde t>i_n}\frac{|x_{\widetilde t}|^p+\cdots+|x_{i_n+1}|^p}{a_{\widetilde t}(p)+\cdots+a_{i_n+1}(p)}\Big\}.
\end{equation}
Then $I_n=\{t\in\mathbb N\,:\, i_{n-1}<t\leq i_n\}$ is a partition of $\mathbb N$ satisfying 
\begin{align}
\label{eq:ineq1}
\sup_{t\in I_n}\frac{|x_{i_{n-1}+1}|^p+|x_{i_{n-1}}|^p+\cdots+|x_t|^p}{a_{i_{n-1}+1}(p)+a_{i_{n-1}}(p)+\cdots+a_t(p)}\leq \frac{\displaystyle\sum_{k\in I_n}|x_k|^p}{\displaystyle\sum_{k\in I_n}a_k(p)}
\end{align}
and
\begin{align}
\label{eq:ineq2}
\frac{\displaystyle\sum_{k\in I_n}|x_k|^p}{\displaystyle\sum_{k\in I_n}a_k(p)}>\frac{\displaystyle\sum_{k\in I_{n+1}}|x_k|^p}{\displaystyle\sum_{k\in I_{n+1}}a_k(p)}.
\end{align}

\end{lemma}
\begin{proof}
Note that the denominator in \eqref{eq:ratios} is always positive since $a_{i_n+1}$ cannot vanish. It is evident for $n=0$ and for $n=1,2,3,\cdots$ we have the inequality 
$$\frac{\displaystyle\sum_{k\in I_n}|x_k|^p}{\displaystyle\sum_{k\in I_n}a_k(p)}>\frac{\displaystyle\sum_{k\in I_n}|x_k|^p+|x_{i_n+1}|^p}{\displaystyle\sum_{k\in I_n}a_k(p)+a_{i_n+1}(p)}$$
since otherwise $i_n$ could not be the last time at which the ratio on the left side attained its maximum. Moreover $x\in\ell^p$ implies that \eqref{eq:ratios} is always bounded. If the supremum of \eqref{eq:ratios} is attained finitely often then it is clear how we get the term $i_{n+1}$. But it could happen that the supremum is attained infinitely often or never. In both cases we set $i_{n+1}=\infty$ and so we obtain a finite collection of sets $I_n$. It is clear that in any of the cases whether or not the sequence $(i_n)$ terminates or not we have that $I_n$ form a partition of $\mathbb N$ since $I_j\cap I_k=\emptyset$ when $j\neq k$, else they coincide, and $\cup_{n\in\mathbb N}I_n=\mathbb N$. Inequality \eqref{eq:ineq1} follows immediately from \eqref{eq:ratios}. As for inequality \eqref{eq:ineq2} note that it is equivalent to
$$\frac{\displaystyle\sum_{k\in I_n}|x_k|^p}{\displaystyle\sum_{k\in I_n}a_k(p)}>\frac{\displaystyle\sum_{k\in I_n\cup I_{n+1}}|x_k|^p}{\displaystyle\sum_{k\in I_n\cup I_{n+1}}a_k(p)}$$
that itself follows from \eqref{eq:ratios}. This completes the proof. 
\end{proof}

\subsection{Factorization of sequence spaces}
This part is an excursion into the factorization theorems of Bennett \cite{Bennett}. With minor modifications his results and techniques apply to spaces $\ell^p_M$ for a large class of infinite matrices $M$. 
We investigate the problem of identification of $(\ell^p_M)^*$ by means of factorization of sequence spaces. 
To this end let $p\geq 1$ and $M$ be an infinite matrix with $\ell^p$-summable diagonal terms. We associate to $M$ a nonnegative sequence of numbers
\begin{equation}
\label{eq:a-sequence}
a_n(p)=\Big(\sum_{k\geq n}|m_{kk}|^p\Big)^{-1}-\Big(\sum_{k\geq n-1}|m_{kk}|^p\Big)^{-1},\quad n\geq 2
\end{equation}
with
\begin{equation}
\label{eq:a-sequence2}
a_1(p)=\Big(\sum_{k=1}^{\infty}|m_{kk}|^p\Big)^{-1}.
\end{equation}
By construction we have 
\begin{equation}
\label{eq:a-sequence3}
\Big(\sum_{k\geq n}|m_{kk}|^p\Big)^{-1}=\sum_{j=1}^{n}a_j(p),\quad n\in\mathbb N.
\end{equation}
Denote by $A_n(p)=\sum_{j=1}^{n}a_j(p)$.
In the spirit of \cite[Bennett]{Bennett} we introduce the spaces
\begin{equation}
\label{eq:d(p)}
d_M(p)=\{x\in\mathbb C^{\mathbb N}\,:\, \|x\|_{d_M(p)}\coloneqq\Big(\sum_{n=1}^{\infty}a_n(p)\,\sup_{k\geq n}|x_k|^p\Big)^{1/p}<\infty\}.
\end{equation}
and 
\begin{equation}
\label{eq:g(q)}
g_M(q)=\{x\in\mathbb C^{\mathbb N}\,:\,\|x\|_{g_M(q)}:=\sup_{n}A^{-1/p}_n(p)\Big(\sum_{k=1}^n|x_k|^q\Big)^{1/q}<\infty\}.
\end{equation}
By similar arguments as in Theorem \ref{th:M-Banach-space} it can be shown that $(d_M(p),\|\cdot\|_{d_M(p)})$ and $(g_M(q),\|\cdot\|_{g_M(q)})$ are Banach spaces. It is not required at this point that $p^{-1}+q^{-1}=1$, but only that $p,q\geq 1$. In particular in the definition of $g_M(q)$ we always normalize the inner expression in \eqref{eq:g(q)} by the reciprocal of $A_n^{1/p}(p)$, the $\ell^p-$norm of the tail of the diagonal terms of $M$. Note that while in \cite[\S3, \S12]{Bennett} the sequence $(a_n)$ is arbitrary, here it is determined by the underlying matrix $M$, therefore making a natural connection between $M$ and its associated sequence space $\ell^p_M$.
Given two sequences $y=(y_n), z=(z_n)$ let $y\cdot z=(y_n\cdot z_n)$.

\begin{thm}\cite[Theorem 3.8]{Bennett}
\label{th:d-g}
Let $0<p\leq \infty$, then $d_M(p)\cdot g_M(p)=\ell^p$.  In particular the relation holds
\begin{equation}
\label{eq:inf-relation}
\|x\|_p= \inf\{\|z\|_{g_M(p)}\,\|y\|_{d_M(p)}\,:\,y\cdot z=x\}
\end{equation}
\end{thm}
\begin{proof}
Let $y\in d_M(p)$ and $z\in g_M(p)$. Denote by $\widehat y_n=\sup_{k\geq n}|y_k|$ the least decreasing majorant of $y$. Let $x=y\cdot z$, then 
\begin{align*}
\|x\|^p_p=\sum_{n=1}^{\infty}|y_n|^p\,|z_n|^p\leq \sum_{n=1}^{\infty}\widehat y^p_n\,|z|^p.
\end{align*}
On the other hand we have
\begin{align*}
\sum_{k=1}^n|z_k|^p\leq \|z\|_{g_M(p)}^p\,A_n(p)=\|z\|_{g_M(p)}^p\,\sum_{k=1}^na_k(p).
\end{align*}
Then by Lemma \ref{l:partial-sums} we obtain 
\begin{align*}
\sum_{k=1}^n\widehat y^p_k\,|z_k|^p\leq \|z\|_{g_M(p)}^p\,\sum_{k=1}^n\widehat y_k^p\,a_k(p)
\end{align*}
consequently
\begin{align*}
\|x\|^p_p\leq \sum_{n=1}^{\infty}\widehat y_n^p\,|z_n|^p\leq \|z\|_{g_M(p)}^p\,\sum_{n=1}^{\infty}\widehat y_n^p\,a_n(p)=\|z\|_{g_M(p)}^p\,\|y\|^p_{d_M(p)}<\infty.
\end{align*}
This proves $d_M(p)\cdot g_M(p)\subseteq \ell^p$ and thus $\|x\|_p\leq \inf\{\|z\|_{g_M(p)}\,\|y\|_{d_M(p)}\,:\,y\cdot z=x\}$. For the other direction define $y\in\mathbb C^{\mathbb N}$ as follows
$$y_j=\Big(\frac{\displaystyle \sum_{k\in I_n}|x_k|^p}{\displaystyle\sum_{k\in I_n}a_k(p)}\Big)^{1/p},\quad j\in I_n.$$
Then by Lemma \ref{l:Bennett}, \eqref{eq:ineq2} we have that $(y_j)$ is a nonincreasing sequence of nonnegative numbers. Consequently 
\begin{align*}
\|y\|^p_{d_M(p)}=\sum_{n=1}^{\infty}a_n(p)\,\widehat y^p_n=\sum_{n=1}^{\infty}a_n(p)\, y^p_n=\sum_{n=1}^{\infty}\sum_{k\in I_n}a_k(p)\,y_k^p=\sum_{n=1}^{\infty}\sum_{k\in I_n}|x_k|^p=\|x\|^p_{p}.
\end{align*}
Next define $z\in\mathbb C^{\mathbb N}$ by
\begin{align*}
z_j=\Big(\frac{\displaystyle\sum_{k\in I_n}a_k(p)}{\displaystyle \sum_{k\in I_n}|x_k|^p}\Big)^{1/p}\,x_j,\quad j\in I_n.
\end{align*}
Evidently it holds that $x=y\cdot z$ and moreover when $k\in I_n$ we get 
\begin{align*}
\sum_{j=1}^k|z_j|^p&=\sum_{j\in I_1\cup\cdots\cup I_{n-1}}|z_j|^p+\Big(|z_{i_{n-1}+1}|^p+\cdots+|z_k|^p\Big)\\&
=\sum_{j\in I_1\cup\cdots\cup I_{n-1}}a_j(p)+\Big(|x_{i_{n-1}+1}|^p+\cdots+|x_k|^p\Big)\,\frac{\displaystyle\sum_{j\in I_n}a_j(p)}{\displaystyle \sum_{j\in I_n}|x_j|^p}
\leq \sum_{j=1}^ka_j(p).
\end{align*}
In view of formula \eqref{eq:a-sequence3} we then have
\begin{align*}
A_k^{-1}(p)\,\sum_{j=1}^k|z_k|^p\leq 1
\end{align*}
i.e. $\|z\|_{g_M(p)}\leq 1$. So $\ell^p\subseteq d_M(p)\cdot g_M(p)$ and $ \inf\{\|z\|_{g_M(p)}\,\|y\|_{d_M(p)}\,:\,y\cdot z=x\}\leq \|x\|_p$. 
\end{proof}

\begin{thm}\cite[Theorem 12.3]{Bennett}
\label{th:dual}
Let $p,q>1$ satisfy $p^{-1}+q^{-1}=1$ , then the following relations hold
\begin{equation}
\label{eq:relation1}
d^*_M(p)=\ell^q\cdot g_M(p)
\end{equation}
and
\begin{equation}
\label{eq:relation2}
(\ell^q\cdot g_M(p))^*=d_M(p).
\end{equation}
\end{thm}
\begin{proof}
First we show \eqref{eq:relation1}. Let $y\in\ell^q,\,z\in g_M(p)$. Denote $x=y\cdot z$ then for any $u\in d_M(p)$ we have 
\begin{align*}
\|x\cdot u\|_1=\|y\cdot z\cdot u\|_1\leq \|y\|_q\,\|z\cdot u\|_p.
\end{align*}
In view of Theorem \ref{th:d-g} we obtain
\begin{align*}
\|x\cdot u\|_1\leq \|y\|_q\,\|z\|_{g_M(p)}\,\|u\|_{d_M(p)}.
\end{align*}
Let $\varphi_u$ be the associated linear functional to $u$, i.e. $\varphi_u(x)=\sum_{n=1}^{\infty}x_n\cdot u_n$, then
\begin{align*}
\sup_{u\neq 0}\frac{|\varphi_u(x)|}{\|u\|_{d_M(p)}}\leq \sup_{u\neq 0}\frac{\|x\cdot u\|_1}{\|u\|_{d_M(p)}}\leq \|y\|_q\,\|z\|_{g_M(p)}
\end{align*}
consequently $x\in(d_M(p))^*$ yielding the inclusion $\ell^q\cdot g_M(p)\subseteq (d_M(p))^*$. Now let $x\in (d_M(p))^*$ and define 
\begin{equation}
\label{eq:psi}
\psi(x)\coloneqq \Big(\sum_{n=1}^{\infty}\Big(\sum_{j\in I_n}a_j\Big)^{1-q}\,\Big(\sum_{j\in I_n}|x_j|\Big)^{q}\Big)^{1/q}
\end{equation}
where $I_n$ is a partition of $\mathbb N$ similarly as constructed in Lemma \ref{l:Bennett}. Note that 
$$\|x\|_{(d_M(p))^*}\geq \sup_{n\in\mathbb N} \Big(\sum_{k=1}^na_k(p)\Big)^{-1/p}\,\Big(\sum_{k=1}^n|x_k|\Big)$$
since $\|x\|_{(d_M(p))^*}=\sup\{\sum_{k=1}^{\infty}|x_k\,y_k|\,:\,\|y\|_{d_M(p)}\leq 1\}\geq \sup_{n\in\mathbb N}\sum_{k=1}^{n}|x_k|\,|y_{k,n}|$, where $y_{k,n}=\Big(\sum_{k=1}^na_k(p)\Big)^{-1/p}$ for $k\leq n$ and $y_{k,n}=0$ else. In particular this means that the partition $I_n$ as in Lemma \ref{l:Bennett} is well-defined, whenever $x\in (d_M(p))^*$.
We claim that $\psi(x)\leq \|x\|_{(d_M(p))^*}$. Consider the sequence 
\begin{align*}
u_k\coloneqq \Big(\sum_{j\in I_n}a_j\Big)^{1-q}\,\Big(\sum_{j\in I_n}|x_j|\Big)^{q-1}\,\quad k\in I_n.
\end{align*}
Then
\begin{align*}
\|x\cdot u\|_1=\sum_{n=1}^{\infty}\sum_{k\in I_n}|x_k|\,|u_k|=\sum_{n=1}^{\infty}\Big(\sum_{j\in I_n}a_j\Big)^{1-q}\,\Big(\sum_{j\in I_n}|x_j|\Big)^{q}=\psi^q(x).
\end{align*}
In view of Lemma \ref{l:Bennett} we have that $(u_k)$ is a nonincreasing sequence which then yields
\begin{align*}
\|u\|^p_{d_M(p)}=\sum_{n=1}^{\infty}\sum_{k\in I_n}a_k\,\sup_{j\geq k}u^p_j=\sum_{n=1}^{\infty}\sum_{k\in I_n}a_k\,u^p_k=\sum_{n=1}^{\infty}\Big(\sum_{j\in I_n}a_j\Big)^{1-q}\,\Big(\sum_{j\in I_n}|x_j|\Big)^{q}=\psi^q(x).
\end{align*}
Consequently
\begin{align*}
\|x\|_{(d_M(p))^*}\geq \sup_{u\neq 0}\frac{\|x\cdot u\|_1}{\|u\|_{d_M(p)}}=\psi(x).
\end{align*}
Next we take $y,\,z\in\mathbb C^{\mathbb N}$ defined by
\begin{align*}
z_j=\Big(|x_j|\,\Big(\sum_{k\in I_n}a_k\Big)\,\Big(\sum_{k\in I_n}|x_k|\Big)^{-1}\Big)^{1/p},\quad j\in I_n
\end{align*}
and $y_j=x_j/z_j$ for $j\in I_n$.
Then clearly $x=y\cdot z$. By similar arguments as in Theorem \ref{th:d-g} we have $\|z\|_{g_{M}(p)}\leq 1$ and so in particular $z\in g_M(p)$. On the other hand we have 
\begin{align*}
\|y\|_q^q=\sum_{n=1}^{\infty}\sum_{k\in I_n}\Big|\frac{x_k}{z_k}\Big|^q&=\sum_{n=1}^{\infty}\sum_{k\in I_n}|x_k|\Big(\sum_{j\in I_n}a_j\Big)^{1-q}\,\Big(\sum_{j\in I_n}|x_j|\Big)^{q-1}\\&=\sum_{n=1}^{\infty}\Big(\sum_{j\in I_n}a_j\Big)^{1-q}\,\Big(\sum_{j\in I_n}|x_j|\Big)^{q}=\psi^q(x)
\end{align*}
implying $\|y\|_q\,\|z\|_{g_M(p)}\leq \psi(x)$, consequently $\|y\|_q\,\|z\|_{g_M(p)}\leq\|x\|_{(d_M(p))^*}$. This proves the inclusion $(d_M(p))^*\subseteq\ell^q\cdot g_M(p)$ and thus $(d_M(p))^*=\ell^q\cdot g_M(p)$. 

Next we demonstrate \eqref{eq:relation2}. First note that $d_M(p)\subseteq (d_M(p))^{**}=(\ell^q\cdot g_M(p))^*$ from relation \eqref{eq:relation1}. Now let $x\in (\ell^q\cdot g_M(p))^*$, then $x\in\ell^{\infty}$ since $\ell^p\subseteq g_M(p)$. Indeed let $n\in\mathbb N$ and let $y=z$ with $y_k=0$ for all $k\neq n$ and $y_k=1$ when $k=n$. Then $\|y\|_{p}=\|y\|_q=1$ implies that $y\in\ell^p\cap\ell^q$ and $\|x\cdot(y\cdot z)\|_1=|x_n|\leq C$ for some $C>0$. 
Next define the sequence of nonnegative integers $(i_n)$ with $i_0=0$ and $i_{n+1}\coloneqq \sup\{t>i_n\,:\,|x_t|=\sup_{k>i_n}|x_k|\}$. Similar to method in Lemma \ref{l:Bennett} the sequence $(i_n)$ is well-defined and the sets $I_n=\{t\in\mathbb N\,:\,i_{n-1}<t\leq i_n\}$ determine a partition of $\mathbb N$. Let $w\in\mathbb C^{\mathbb N}$ be given by
\begin{align*}
w_k=\left\{
\begin{array}{ll}
      \displaystyle\Big(\sum_{j\in I_n}^ka_j(p)\Big)^{1/p} & k= i_n\;\text{for some}\;n \\[1em]
     0 & \text{else}.
\end{array} 
\right.
\end{align*}
Then it can be shown that $w\in g_M(p)$ with $\|w\|_{g_M(p)}=1$. On the other hand we have
\begin{align*}
\|x\|_{(\ell^q\cdot g_M(p))^*}&=\sup\{\|x\cdot y\|_1\,:\,y=u\cdot v,\,\|u\|_q\leq 1,\,\|v\|_{g_M(p)}\leq 1 \}\\&
=\sup\{\|x\cdot v\|_p\,:\,\|v\|_{g_M(p)}\leq 1\}\\&\geq \|x\cdot w\|_p\\&=\Big(\sum_{n=1}^{\infty}\sum_{k\in I_n}|x_k|^p\,|w_k|^p\Big)^{1/p}\\&=\Big(\sum_{n=1}^{\infty}|x_{i_n}|^p\,\sum_{j\in I_n}a_j(p)\Big)^{1/p}\geq \Big(\sum_{n=1}^{\infty}\sum_{j\in I_n}a_j(p)\,\sup_{k\geq j}|x_k|^p\Big)^{1/p}=\|x\|_{d_M(p)}.
\end{align*}
This shows the inclusion $(\ell^q\cdot g_M(p))^*\subseteq d_M(p)$ and \eqref{eq:relation2} is proved. 
\end{proof}

\begin{corollary}
\label{c:factorization-dual}
 $(\ell^p_M)^*$ is isomorphic to $d_M(q)$ whenever the sequence space $\ell^p_M$ satisfies the factorization $\ell^p_M=\ell^p\cdot g_M(q)$. 
\end{corollary}

\subsection{Identification of the dual} In view of Corollary \ref{c:factorization-dual} the question of identification of $(\ell^p_M)^*$ reduces to finding conditions, desirably necessary and sufficient, for which the factorization $\ell^p_M=\ell^p\cdot g_M(q)$ holds. Here we assume that $p,q\geq 1$ and $p^{-1}+q^{-1}=1$.

\begin{thm}
\label{th:factorization}
Let $M=(m_{nk})$ be a lower triangular matrix with $\ell^p$-summable diagonal terms such that $|m_{n(k+1)}|\leq |m_{nk}|$ for all $k\in\{1,2,\cdots,n-1\}$ and all $n\in\mathbb N$.
 Then every $x\in\ell^p_M$ admits a factorization $x=y\cdot z$ with $y\in\ell^p$ and $z\in g_M(q)$.
\end{thm}
\begin{proof}
We follow similar footsteps as in \cite[Theorem 4.5]{Bennett} where it was shown for the case when $M=C$ (the Ces\'aro method). First assume that $p>1$. Let $x\in\ell^p_M$ and let 
\begin{equation}
\label{eq:bn}
b_n\coloneqq \sum_{k=n}^{\infty}|m_{kk}|\,(\sum_{j=1}^k|m_{kj}|\,|x_j|)^{p-1}\quad (n=1,2,3,\cdots).
\end{equation}
Note that $b_n$ is finite for every $n\in\mathbb N$. Indeed by H\"older's inequality and $\ell^p$-summability of the diagonal terms we get
\begin{align*}
b_n\leq\sum_{k=1}^{\infty}|m_{kk}|\,(\sum_{j=1}^k|m_{kj}|\,|x_j|)^{p-1}\leq (\sum_{k=1}^{\infty}|m_{kk}|^p)^{1/p}\,\|x\|_{M}^{p/q}<\infty.
\end{align*}
Also $(b_n)$ is nonincreasing and nonnegative. Let $y_n=x_n\,|x_n|^{1/p-1}\,b_n^{1/p}$ and $z_n=|x_n|^{1/q}\,b_n^{-1/p}$ for $n\in\mathbb N$, then $x=y\cdot z$. Note that
\begin{align*}
\sum_{n=1}^{\infty}|y_n|^p=\sum_{n=1}^{\infty}|x_n|\,b_n&=\sum_{n=1}^{\infty}|x_n|\,\sum_{k=n}^{\infty}|m_{kk}|\,(\sum_{j=1}^k|m_{kj}|\,|x_j|)^{p-1}\\&
=\sum_{n=1}^{\infty}|m_{nn}|\,(\sum_{k=1}^n|x_{k}|)\,(\sum_{k=1}^n|m_{nk}|\,|x_k|)^{p-1}\\&
\leq \sum_{n=1}^{\infty}(\sum_{j=1}^n|m_{nk}|\,|x_k|)^{p}
\end{align*}
where in the last step we have used the monotonicity of $(|m_{nk}|)$ for each $n\in\mathbb N$. Hence from inequality above we arrive at $\|y\|_p\leq \|x\|_{M,p}$ implying $y\in\ell^p$. Next we have
\begin{align*}
(\sum_{k=1}^N|z|^q)^p=(\sum_{k=1}^N|x_k|^{1/q}\,|x_k|^{1/p}\,b_k^{-q/p})^p\leq (\sum_{k=1}^N|x_k|)^{p-1}\,(\sum_{k=1}^{N}|x_k|\,b_k^{-q}).
\end{align*}
For $n\geq N$ we then obtain 
\begin{align*}
\sum_{n=N}^{\infty}|m_{nn}|^p\,(\sum_{k=1}^N|z|^q)^p&\leq \sum_{n=N}^{\infty}|m_{nn}|^p\,(\sum_{k=1}^n|x_k|)^{p-1}\,(\sum_{k=1}^{N}|x_k|\,b_k^{-q})\\&\leq\sum_{n=N}^{\infty}|m_{nn}|\,(\sum_{k=1}^n|m_{nk}|\,|x_k|)^{p-1}\,(\sum_{k=1}^{N}|x_k|\,b_k^{-q})\\&
=b_N\,(\sum_{k=1}^{N}|x_k|\,b_k^{-q})\leq \sum_{k=1}^{N}|x_k|\,b_k^{1-q}=\sum_{k=1}^{N}|z_k|^q
\end{align*}
where we have used monotonicity of $(|m_{nk}|)$ and of $(b_n)$ in the second and third inequality respectively. Consequently 
\begin{align*}
A_N^{-1}(p)\,(\sum_{k=1}^N|z_k|^q)^{p-1}\leq 1
\end{align*}
implying $\|z\|_{g_M(q)}\leq 1$. The case $p=1$ is dealt accordingly by taking $b_n=\sum_{k=n}^{\infty}|m_{kk}|$ and $z_n=b^{-1}_n$ for $n\in\mathbb N$. This completes the proof.
\end{proof}

 To $M$ we associate another sequence of nonnegative numbers
\begin{equation}
\label{eq:b-sequence}
b_n(p,q)=\Big(\sum_{k\geq n}|m_{kk}|^p\Big)^{-q/p}-\Big(\sum_{k\geq n-1}|m_{kk}|^p\Big)^{-q/p},\quad n\geq 2
\end{equation}
with
\begin{equation}
\label{eq:b-sequence2}
b_1(p,q)=\Big(\sum_{k=1}^{\infty}|m_{kk}|^p\Big)^{-q/p}.
\end{equation}
By construction we have 
\begin{equation}
\label{eq:b-sequence3}
\Big(\sum_{k\geq n}|m_{kk}|^p\Big)^{-q/p}=\sum_{j=1}^{n}b_j(p,q),\quad n\in\mathbb N.
\end{equation}
Denote by $B_n(p,q)=\sum_{j=1}^{n}b_j(p,q)$ and $\widehat b_n(p,q)=\sup_{k\leq n}b_k(p,q)$. Evidently the relation $B_n(p,q)=A_n^{q/p}(p)$ holds for every $n\in\mathbb N$. 

\begin{thm}
\label{th:converse}
 Let $M=(m_{nk})$ be a matrix as in the previous Theorem. If additionally 
\begin{equation}
\label{eq:assumption-growth}
m_{n1}= \left\{
        \begin{array}{ll}
           O((n\,\widehat b^{1/q}_n)^{-1}) & p>1 \\
            O((n^{1+\epsilon}\,A_n)^{-1}) & p=1,\;\text{for some}\;\epsilon>0
        \end{array}
    \right.
\end{equation}
then $x\in\ell^p_M$ if and only if $x$ admits a factorization $x=y\cdot z$ where $y\in\ell^p$ and $z\in g_M(q)$. In particular $\ell^p_M=\ell^p\cdot g_M(q)$.
\end{thm}
\begin{proof}
The inclusion $\ell^p_M\subseteq \ell^p\cdot g_M(q)$ follows from Theorem \ref{th:factorization}. We show the reverse inclusion. First suppose that $p>1$. Let $x=y\cdot z$ with $y\in\ell^p$ and $z\in g_M(q)$. 
Notice that 
\begin{align*}
\sum_{k=1}^N|z_k|^q\leq \|z\|_{g_M(q)}^q\,B_N(p,q),\quad (N=1,2,3,\cdots).
\end{align*}
For $p>1$ consider the sequence 
\begin{equation*}
\label{eq:w}
w_k={k-1-1/p\choose k-1},\;k\in\mathbb N.
\end{equation*}
$(w_k)$ is positive, decreasing and satisfies the inequality (see \cite[Lemma 4.11]{Bennett})
\begin{equation}
\label{eq:ineq-w}
(w_1+w_2+\cdots+w_k)^{p-1}<(kq)^p\,(w_k^{p-1}-w_{k+1}^{p-1}),\quad k\in\mathbb N.
\end{equation}
Let $\widetilde b_{nk}(p,q)=b_{\sigma_n(k)}(p,q)$ be the nondecreasing reordering of the set $\{b_1(p,q),\cdots, b_n(p,q)\}$, i.e. $\widetilde b_{n1}(p,q)\leq \widetilde b_{n2}(p,q)\leq\cdots\leq \widetilde b_{nn}(p,q)$. Moreover note that $$\sum_{k=1}^nb_k(p,q)=\sum_{k=1}^n\widetilde b_{nk}(p,q),\quad (n=1,2,3,\cdots).$$
By Lemma \ref{l:partial-sums}, since also the sequence $w_k\,\widetilde b_{nk}(p,q)^{-1}$ is nonincreasing, we have 
\begin{align}
\label{eq:ineq-zk}
\sum_{k=1}^N|z_k|^q\,\frac{w_k}{\widetilde b_{nk}(p,q)}\leq \|z\|_{g_M(q)}^q\,\sum_{k=1}^Nw_k,\quad (N=1,2,3,\cdots).
\end{align}
Then we consider the sum 
\begin{align*}
\Big(\sum_{k=1}^{n}|m_{nk}|\,|x_k|\Big)^p&=\Big(\sum_{k=1}^{n}|m_{nk}|\,|y_k|\, \frac{\widetilde b^{1/q}_{nk}(p,q)}{w_k^{1/q}}\,|z_k|\, \frac{w_k^{1/q}}{\widetilde b^{1/q}_{nk}(p,q)}\Big)^p\\&\leq \Big(\sum_{k=1}^n|m_{nk}|^p\,|y_k|^p\, \frac{\widetilde b^{p/q}_{nk}(p,q)}{w_k^{p/q}}\Big)\,\Big(\sum_{k=1}^n|z_k|^q\,\frac{w_k}{\widetilde b_{nk}(p,q)}\Big)^{p/q}.
\end{align*}
By assumption $|m_{nk}|\leq |m_{n1}|$ for $k\leq n$. Using inequalities \eqref{eq:ineq-w} and \eqref{eq:ineq-zk} we then obtain 
\begin{align*}
\Big(\sum_{k=1}^{n}|m_{nk}|\,|x_k|\Big)^p&\leq |m_{n1}|^p\,\|z\|_{g_M(q)}^p\,\Big(\sum_{k=1}^n|y_k|^p\, \frac{\widetilde b^{p/q}_{nk}(p,q)}{w_k^{p/q}}\Big)\,\Big(\sum_{k=1}^nw_k\Big)^{p-1}\\&
\leq |m_{n1}|^p\,\|z\|_{g_M(q)}^p\,\Big(\sum_{k=1}^n|y_k|^p\, \frac{\widetilde b^{p/q}_{nk}(p,q)}{w_k^{p/q}}\Big)\,\Big(\sum_{k=1}^nw_k\Big)^{p-1}
\\&
\leq (n\,q)^p\,|m_{n1}|^p\,\|z\|_{g_M(q)}^p\,\Big(\sum_{k=1}^n|y_k|^p\, \frac{\widetilde b^{p/q}_{nk}(p,q)}{w_k^{p/q}}\Big)\,(w_n^{p-1}-w_{n+1}^{p-1}).
\end{align*}
Realizing that $\widetilde b_{nn}(p,q)=\sup_{k\leq n}b_k(p,q)=\widehat b_n(p,q)$ and noting that $p/q=p-1$ we obtain 
\begin{align*}
\sum_{n=1}^{\infty}\Big(\sum_{k=1}^{n}|m_{nk}|\,|x_k|\Big)^p&\leq q^p\,\|z\|_{g_M(q)}^p\,\sum_{n=1}^{\infty}(n\,|m_{n1}|\,\widehat b^{1/q}_n(p,q))^p\,(w_n^{p-1}-w_{n+1}^{p-1})\,\sum_{k=1}^n|y_k|^p\, w_k^{-(p-1)}\\&
=q^p\,\|z\|_{g_M(q)}^p\,\sum_{n=1}^{\infty}|y_n|^p\, w_n^{-(p-1)}\,\sum_{k=n}^{\infty}(k\,|m_{k1}|\,\widehat b^{1/q}_k(p,q))^p\,(w_k^{p-1}-w_{k+1}^{p-1}).
\end{align*}
By assumption \eqref{eq:assumption-growth} there is $N\in\mathbb N$ large enough so that $|m_{k1}|\leq C\,k^{-1}\,\widehat b^{-1}_k(p,q)$ for all $k\geq N$ and some $C>0$. We obtain the following upper estimate
\begin{align*}
\sum_{n\geq N}|y_n|^p\, w_n^{-(p-1)}\,\sum_{k=n}^{\infty}(k\,|m_{k1}|\,\widehat b^{1/q}_k(p,q))^p&\,(w_k^{p-1}-w_{k+1}^{p-1})\\&\leq C^p\,\sum_{n\geq N}|y_n|^p\, w_n^{-(p-1)}\,\sum_{k=n}^{\infty}(w_k^{p-1}-w_{k+1}^{p-1})\\&
=C^p\,\sum_{n\geq N}|y_n|^p<\infty.
\end{align*}
For the rest of the sum we have the upper estimate
\begin{align*}
\sum_{n=1}^{N-1}|y_n|^p\, w_n^{-p/q}\,&\sum_{k=n}^{\infty}(k\,|m_{k1}|\,\widehat b^{1/q}_k(p,q))^p\,(w_k^{p-1}-w_{k+1}^{p-1})\\&=\sum_{n=1}^{N-1}|y_n|^p\, w_n^{-p/q}\,\sum_{k=n}^{N-1}(k\,|m_{k1}|\,\widehat b^{1/q}_k(p,q))^p\,(w_k^{p-1}-w_{k+1}^{p-1})\\&
+\sum_{n=1}^{N-1}|y_n|^p\, w_n^{-p/q}\,\sum_{k\geq N}(k\,|m_{k1}|\,\widehat b^{1/q}_k(p,q))^p\,(w_k^{p-1}-w_{k+1}^{p-1}).
\end{align*}
The first sum is finite as only finitely many terms are present. The second sum can be estimated from above by
\begin{align*}
\leq C^p\,\sum_{n=1}^{N-1}|y_n|^p\, w_n^{-(p-1)}\,\sum_{k\geq N}(w_k^{p-1}-w_{k+1}^{p-1})=C^p\,\sum_{n=1}^{N-1}|y_n|^p\frac{w_N^{p-1}}{w_n^{p-1}}\leq C^p\,\sum_{n=1}^{N-1}|y_n|^p.
\end{align*}
Therefore all in all we obtain that 
$$\sum_{n=1}^{\infty}\Big(\sum_{k=1}^{n}|m_{nk}|\,|x_k|\Big)^p<\infty.$$
Last we consider the special case $p=1$. Take $x=y\cdot z$ with $z\in\ell^1$ and $z\in g_M(\infty)$, i.e. $\sup_{n\in\mathbb N}(A_n^{-1}(1)\,\sup_{k\leq n}|z_k|)<\infty$. It follows that $\sup_{k\leq N}|z_k|\leq \|z\|_{g_M(\infty)}\,A_N(1)$ for all $N=1, 2, 3, \cdots$. We then have the estimate
\begin{align*}
\sum_{n=1}^{\infty}\sum_{k=1}^{n}|m_{nk}|\,|x_k|&=\sum_{n=1}^{\infty}\sum_{k=1}^{n}|m_{nk}|\,|y_k|\,|z_k|\leq \sum_{n=1}^{\infty}\sup_{k\leq n}|z_k|\,\sum_{k=1}^{n}|m_{nk}|\,|y_k|\\&\leq \|z\|_{g_M(\infty)}\, \sum_{n=1}^{\infty}A_n(1)\,\sum_{k=1}^{n}|m_{nk}|\,|y_k|
\\&\leq \|z\|_{g_M(\infty)}\, \|y\|_1\sum_{n=1}^{\infty}A_n(1)\,|m_{n1}|\\&
\leq C_M\, \|z\|_{g_M(\infty)}\, \|y\|_1\sum_{n=1}^{\infty}\frac{1}{n^{1+\epsilon}}=C_M\, \|z\|_{g_M(\infty)}\, \|y\|_1\,\zeta(1+\epsilon)<\infty
\end{align*}
as $\epsilon>0$. The positive constant $C_M$ depends on $M$. This completes the proof.
\end{proof}

\begin{corollary}
\label{c:dual-identification}
Under the conditions in Theorem \ref{th:converse} and in view of Corollary \ref{c:factorization-dual} for all $p> 1$ the dual $(\ell^p_M)^*$ is isometric isomorphic to $d_M(q)$.
\end{corollary}

\begin{example}
\label{ex:cesaro-alpha}
Consider the class of (generalized) Ces\'aro matrix $C_{\alpha}=(c_{nk})$ for $0<\alpha\leq 1$ given by \eqref{eq:Cesaro}.
Then 
\begin{align*}
c_{nk}&=\frac{(n-k+\alpha-1)\cdots(\alpha+1)\,\alpha}{(n-k)!}\cdot\frac{(n-1)!}{(n+\alpha-1)\cdots(\alpha+1)}\\&
=\frac{n-1}{n+\alpha-2}\,\frac{n-2}{n+\alpha-3}\cdots\frac{n-k+1}{n+\alpha-k}\,\frac{\alpha}{n+\alpha-1}
\end{align*}
implies in particular $c_{n(k+1)}\leq c_{nk}$ for all $k\in\{1,2,\cdots,n-1\}$. Moreover by Stirling's approximation we have $$c_{nn}=\frac{(n-1)!}{(n+\alpha-1)\cdots(\alpha+1)}\sim \text{const.}\, \frac{1}{n^{\alpha}},$$
therefore $\sum_{n=1}^{\infty}c_{nn}^p\sim\text{const.}\, \zeta(\alpha\,p)<\infty$ whenever $p>1/\alpha$. For large enough $n\in\mathbb N$ we have 
\begin{align*}
\sum_{k\geq n}c_{kk}^p\sim \sum_{k\geq n}\frac{1}{k^{\alpha\,p}}>\frac{1}{\alpha\,p-1}\,\frac{1}{n^{\alpha\,p-1}}
\end{align*}
and
\begin{align*}
\sum_{k\geq n}c_{kk}^p\sim \sum_{k\geq n}\frac{1}{k^{\alpha\,p}}<\frac{1}{\alpha\,p-1}\,\frac{1}{(n-1)^{\alpha\,p-1}}
\end{align*}
implying 
\begin{align*}
b_n(p,q)<(\alpha\,p-1)^{q/p}(n^{(\alpha-1)q+1}-(n-1)^{(\alpha-1)\,q+1})\sim \text{const.}\,n^{(\alpha-1)\,q}.
\end{align*}
 On the other hand $c_{n1}=\alpha/(n+\alpha-1)$ for every $n\in\mathbb N$ and $\alpha\leq 1$ yield
$c_{n1}=O((n\,\widehat b^{1/q}_n)^{-1})$ for sufficiently large $n\in\mathbb N$. Therefore $\ell^p_{C_{\alpha}}=\ell^p\cdot d_{C_{\alpha}}(q)$ whenever $p>1/\alpha$ and $0<\alpha\leq 1$.
\end{example}

\begin{example}
\label{eq:power-type}
Let $p>1$ and $M=(m_{nk})$ a lower triangular matrix with $|m_{nk}|=\gamma\,n^{-\beta}$ for all $k\leq n$ and every $n\in\mathbb N$ for some $\beta>1/p,\gamma>0$. Obviously $|m_{n(k+1)}|\leq |m_{nk}|$ for all $k\in\{1,2,\cdots, n-1\}$ and $\sum_{n=1}^{\infty}|m_{nn}|^p=\gamma^{-1}\,\sum_{n=1}^{\infty}n^{\beta\,p}=\gamma^{-1}\,\zeta(\beta\,p)<\infty$. By similar estimations as in the previous example we have that $b_n(p,q)\sim\text{const.}\,n^{(\beta-1)\,q}$ implying $|m_{n1}|=\gamma^{-1}\,n^{-\beta}=O((n\,\widehat b^{1/q}_n)^{-1})$. Therefore again conclusion of Theorem \ref{th:converse} holds true.
\end{example}

\begin{remark}
\label{r:Cesaro} 
In both examples the Ces\'aro method  proven by \cite[Bennett]{Bennett} is a special case. In the first example with $\alpha=1$ and in the second with $\beta=\gamma=1$.
\end{remark}

\bibliographystyle{plain}
\bibliography{references}

\end{document}